\documentclass[a4paper,pdftex,reqno,10pt]{amsart}
\usepackage{geometry}
\usepackage{graphicx}
\usepackage{amssymb,amsmath}

\def\Ps{\mathcal{P}}

\newcommand{\Tsm}{\hspace*{0.6cm}}

\newcommand{\ii}{i}

\newcommand{\I}{\mathcal{I}}
\newcommand{\IC}{{\dot{\I}(n,2)}}

\def\@begintheorem#1#2{\list{}{\thm@body}%
  \item[]{\bf #1~#2.}\quad\it\ignorespaces}
\def\@opargbegintheorem#1#2#3{\list{}{\thm@body}%
  \item[]{\bf #1~#2~\ifrembrks #3\global\rembrksfalse\else (#3)\fi.}%
  \quad\it\ignorespaces}
\def\@endtheorem{\endlist}
\newtheorem{theorem}{Theorem}
\newtheorem{algo}{Algorithm}
\newtheorem{corollary}{Corollary}
\newtheorem{lemma}{Lemma}

\newtheorem{definition}{Definition}
\newtheorem{conjecture}{Conjecture}
\newcommand{\eop}{\hfill{$\Box$}}

\begin{document}

\title[Integral point sets over finite fields]{Integral point sets over finite fields}
\author{Sascha Kurz}
\address{Sascha Kurz\\Fakult\"at f\"ur Mathematik, Physik und Informatik\\Universit\"at Bayreuth\\Germany}
\email{sascha.kurz@uni-bayreuth.de}

\begin{abstract}
    We consider point sets in the affine plane $\mathbb{F}_q^2$ where each Euclidean distance of two points is an
    element of $\mathbb{F}_q$. These sets are called integral point sets and were  originally defined in $m$-dimensional
    Euclidean spaces $\mathbb{E}^m$. We determine their maximal cardinality $\mathcal{I}(\mathbb{F}_q,2)$. For arbitrary
    commutative rings $\mathcal{R}$  instead of $\mathbb{F}_q$ or for further restrictions as no three points on a line or
    no four points on a circle we give partial results. Additionally we study the geometric structure of the examples with
    maximum cardinality.
\end{abstract}

\maketitle

\section{Introduction}
\noindent
Originally integral point sets were defined in $m$-dimensional Euclidean spaces $\mathbb{E}^m$ as a set of $n$ points with pairwise integral distances in the Euclidean metric, see \cite{integral_distances_in_point_sets,kreisel,phd_kurz,sascha-alfred} for a overview on the most recent results. Here we transfer the concept of an integral point set to modules $\mathcal{R}^m$ of a commutative ring with $1$. We equip those spaces with a squared distance 
$$
  d^2(u,v):=\sum_{i=1}^m (u_i-v_i)^2\quad\in\mathcal{R}.
$$
for any two points $u=(u_1,\dots,u_m)$, $v=(v_1,\dots,v_m)$ in $\mathcal{R}^m$ and say that they are at integral distance if $d^2(u,v)$ is contained in the set $\square_{\mathcal{R}}:=\{r^2\mid r\in\mathcal{R}\}$ consisting of the squares in $\mathcal{R}$. A set of points $\Ps$ is called an integral point set if every pair of points is at integral distance. 

The concept of integral point sets over finite fields is not brand-new. There are some recent papers and preprints \cite{quadrance,Vinh:math0510092,Vinh:math0509598,Vinh:math0606482} by L.A. Vinh dealing with Quadrance graphs. These are in the authors definition point sets in the affine plane $\mathbb{F}_q^2$ where the squared distances, there called quadrances, are elements of $\square_{\mathbb{F}_q}\backslash\{0\}$. So for $q\equiv 3\mod 4$ quadrance graphs coincide with integral point sets over $\mathbb{F}_q^2$. For $q\equiv 1\mod 4$ we have the small difference that $0=0^2$ is not considered as an integral distance. So i.e. the points $(0,0)$ and $(2,3)$ in $\mathbb{F}_{13}$ are not considered to be at an integral distance since $d^2((0,0),(2,3))=2^2+3^2=0$. We would like to mention that quadrance graphs and so integral point sets over finite fields are isomorphic to strongly regular graphs and that there are some connections to other branches of Combinatorics including Ramsey theory and association schemes \cite{0495.05047,0843.05103,1091.05048}. The origin of quadrance graphs lies in the more general concept of rational trigonometry and universal geometry by N.J. Wildberger, see \cite{05044388} for more background.

Some related results on integral point sets over commutative rings can be found in \cite{algo,Dimiev-Setting,axel_1}.

A somewhat older topic of the literature is also strongly connected to integral point sets over finite fields. The Paley graph $\mathcal{PG}_q$ has the elements of the finite field $\mathbb{F}_q$ as its vertices. Two vertices $u$ and $v$ are connected via an edge if and only if their difference is a non-zero square in $\mathbb{F}_q$. For $q=q'^2$ with $q'\equiv 3\mod 4$ we have a coincidence between the Paley graph $\mathcal{PG}_q$ and integral point sets over $\mathcal{PG}_{q'}^2$ or quadrance graphs. It is somewhat interesting that these one-dimensional and two-dimensional geometrical objects are so strongly connected. See i.e. \cite{0876.05093,Vinh:math0509598} for a detailed description and proof of this connection. Actually one uses the natural embedding of $\mathbb{F}_{q^2}$ in $\mathbb{F}_q^2$.

So what are the interesting questions about integral point sets over finite fields? From the combinatorial point of view one could ask for the maximum cardinality $\mathcal{I}(\mathcal{R},m)$ of those point sets in $\mathcal{R}^m$. For $\mathcal{R}=\mathbb{F}_q$ with $q\equiv 3\mod 4$ and $m=2$ this is a classical question about maximum cliques of Paley graphs of square order, where the complete answer is given in \cite{0561.12009}. See also \cite{0945.51004} for some generalizations. A geometer might ask for the geometric structure of the maximal examples. Clearly the case where $\mathcal{R}$ is a finite field $\mathbb{F}_q$ is the most interesting one. 

\subsection{Our contribution}

For primes $p$ we completely classify maximal integral point sets in the affine planes $\mathbb{F}_p^2$ and for prime powers $q=p^r$ we give partial results. Since in an integral point set not all directions can occur we can apply some R\'edei-type results in this context. Although these results are not at hand in general we can derive some results for arbitrary rings $\mathcal{R}$ and special cases like $\mathcal{R}=\mathbb{Z}_{p^2}$ or rings with characteristic two. 

It will turn out that most maximal examples or constructions in the plane consist of only very few lines. So it is interesting to consider the case where we forbid three points to be collinear. This means that we look at $2$-arcs with the additional integrality condition. Here we denote the maximal cardinality by $\overline{\mathcal{I}}(\mathcal{R},m)$ where we in general forbid that $m+1$ points are contained in a hyperplane. We give a construction and a conjecture for the case $\mathcal{R}=\mathbb{F}_q$, $2\nmid q$, and $m=2$ using point sets on circles.

Being even more restrictive we also forbid $m+2$ points to be situated on a hypersphere and denote the corresponding maximal cardinality by $\dot{\mathcal{I}}(\mathcal{R},m)$. Although in this case we have almost no theoretical insight so far, this is the most interesting situation when we look from the viewpoint of integral point sets in $\mathbb{E}^m$. As a motivation for further research the following open problem of P.~Erd\H{os} and C.~Noll \cite{0676.52006} may serve:
\begin{quote}
  Are there seven points in the plane, no three on a line, no four on a circle with integral coordinates and pairwise 
  integral distances? 
\end{quote}
If we drop the condition of integral coordinates the problem was recently solved in \cite{kreisel}. As a connection to our problem one may use the ring homomorphism $\mathbb{Z}^m\rightarrow \mathbb{Z}_n^m$, $x\mapsto x+(n\mathbb{Z})^m$, which preserves integral distances and coordinates. For lines and circles the situation is a bit more complicated. We give some examples for various primes $p$ showing $\dot{\mathcal{I}}(\mathbb{Z}_p,2)\ge 7$ and determine some exact numbers. Perhaps in the future an application of the Chinese remainder theorem helps to construct the desired example in $\mathbb{Z}^2$.

\subsection{Organization of the paper}

The paper is arranged as follows. In Section \ref{sec_integral_point_sets} we give the basic definitions and facts on integral point sets over commutative rings $\mathcal{R}$. In Section \ref{sec_automorphism_group} we determine the automorphism group of the affine plane $\mathbb{F}_q^2$ with respect to $\Delta$. For $q\equiv 3\mod 4$ it is the well known automorphism group of the Paley graph $\mathcal{PG}_{q^2}$ which is isomorphic to a subgroup of $\text{PG$\Gamma$}(1,q^2)$ of index $2$, see i.e. \cite{MR0117223,2006math_5252K,1057.05042}. For $q\equiv 1\mod 4$ the automorphism group was not known. We give a proof for both cases and prove some lemmas on integral point sets over finite fields which will be useful in the following sections. Most of the automorphisms also exist in some sense for arbitrary commutative rings $\mathcal{R}$. In Section \ref{sec_maximum_integral_point_sets} we determine the maximum cardinality $\mathcal{I}(\mathbb{F}_q,2)$ of an integral point set over $\mathbb{F}_q^2$ and classify the maximal examples up to isomorphism in some cases. Here we use a result of Blokhuis et al. on point sets with a restricted number of directions. In Section \ref{sec_maximal_integral_point_sets_modular} we give some results on $\mathcal{I}(\mathbb{Z}_n,2)$ and give some constructions which reach this upper bound. In Section \ref{section_arcs} we determine the maximum cardinality $\overline{\mathcal{I}}(\mathbb{F}_q,2)$ of integral point sets over $\mathbb{F}_q$ where no three points are collinear for $q\equiv 3\mod 4$. For $q\equiv 1\mod 4$ we give lower and upper bounds which are only two apart. In Section \ref{sec_general_position} we consider the maximum cardinality $\dot{\mathcal{I}}(\mathbb{F}_q,2)$ of integral point sets over $\mathbb{F}_q^2$ where no three points are collinear and no four points are situated on a circle. We determine some exact values via an exhaustive combinatorial search and list some maximum examples.


\section{Integral point sets}
\label{sec_integral_point_sets}

\noindent
If not stated otherwise we assume that $\mathcal{R}$ is a commutative ring with $1$ and consider sets of elements of 
the $\mathcal{R}$-module $\mathcal{R}^m$. We speak of these elements as points with a geometric interpretation in mind. For our purpose we equip the module $\mathcal{R}^m$ with something similar to an Euclidean metric:

\begin{definition}
  For two points $u=(u_1,\dots,u_m)$, $v=(v_1,\dots,v_m)$ in $\mathcal{R}^m$ we define the \textbf{squared distance} 
  as
  $$
    d^2(u,v):=\sum_{i=1}^m (u_i-v_i)^2\quad\in\mathcal{R}.
  $$
\end{definition}

We are interested in those cases where  $d^2(u,v)$ is contained in the set $\square_{\mathcal{R}}:=\{r^2\mid r\in\mathcal{R}\}$ of squares of $\mathcal{R}$. 

\begin{definition}
  Two points $u=(u_1,\dots,u_m)$, $v=(v_1,\dots,v_m)$ in $\mathcal{R}^m$ are at \textbf{integral distance} if there 
  exists an element $r$ in $\mathcal{R}$ with $d^2(u,v)=r^2$. As a shorthand we define
  $\Delta:\mathcal{R}^m\times\mathcal{R}^m\rightarrow\{0,1\},$
  $$
    (u,v)\mapsto \left\{\begin{array}{ll} 1 & \text{if $u$ and $v$ are at integral distance},\\
    0 & \text{otherwise}.\end{array}\right.
  $$
  A set $\Ps$ of points in $\mathcal{R}^m$ is called an \textbf{integral point set} 
  if all pairs of points are at integral distance. 
\end{definition}

If $\mathcal{R}$ is a finite ring it makes sense to ask for the maximum cardinality of an integral point set in $\mathcal{R}^m$.

\begin{definition}
  By $\mathcal{I}(\mathcal{R},m)$ we denote the maximum cardinality of an integral point set in $\mathcal{R}^m$.
\end{definition}

\begin{lemma}
  $$
    |\mathcal{R}|\le\mathcal{I}(\mathcal{R},m)\le |\mathcal{R}|^m.
  $$
\end{lemma}
\begin{proof}
  For the lower bound we consider the \textit{line} $\Ps=\{(r,0,\dots,0)\mid r\in\mathcal{R}\}$.
\end{proof}

\begin{lemma}
  If $\mathcal{R}$ has characteristic $2$, meaning that $1+1=0$ holds, then we have
  $\mathcal{I}(\mathcal{R},m)=|\mathcal{R}|^m$.
\end{lemma}
\begin{proof}
  For two points $u=(u_1,\dots,u_m)$, $v=(v_1,\dots,v_m)$ in $\mathcal{R}^m$ we have
  $$
    d^2(u,v)=\sum_{i=1}^m (u_i-v_i)^2=\underset{\in\mathcal{R}}{\underbrace{\left(\sum_{i=1}^r u_i+v_i\right)}}^2.
  $$
\end{proof}

So in the remaining part of this article we consider only rings with characteristic not equal to two. If a ring $\mathcal{R}$ is the Cartesian product of two rings $\mathcal{R}_1$, $\mathcal{R}_2$, where we define the operations componentwise, then we have the following theorem:

\begin{theorem}
  \label{thm_decompose}
  $$
    \mathcal{I}(\mathcal{R}_1\times\mathcal{R}_2,m)=\mathcal{I}(\mathcal{R}_1,m)\cdot\mathcal{I}(\mathcal{R}_2,m).
  $$
\end{theorem}
\begin{proof}
  If $\Ps$ is an integral point set in $\mathcal{R}_1\times\mathcal{R}_2$ then the projections into $\mathcal{R}_1$ and 
  $\mathcal{R}_2$  are also integral point sets. If on the other hand $\mathcal{P}_1$ and $\mathcal{P}_2$ are integral 
  point sets over $\mathcal{R}_1$ and $\mathcal{R}_2$, respectively, then $\Ps:=\mathcal{P}_1\times\mathcal{P}_2$ is an 
  integral point set over $\mathcal{R}_1\times\mathcal{R}_2$.
\end{proof}

\begin{lemma}
  If $N$ is an additive subgroup of $\{n\in\mathcal{R}\mid n^2=0\}$ or $\{n\in\mathcal{R}\mid 2n^2=0\,\wedge\,n^2=n^4\}$ 
  then we have for $m\ge 2$ 
  $$
    |N|^{m-1}\cdot|\mathcal{R}|\le\mathcal{I}(\mathcal{R},m)\le|\mathcal{R}|^m.
  $$
\end{lemma}
\begin{proof}
  We can take the integral point set $\Ps=\{(r,n_1,\dots,n_{m-1})\mid r\in \mathcal{R},\,n_i\in N\}$ and 
  have $r^2+\sum\limits_{i=1}^{m-1}n_i^2=r^2$ or
  $r^2+\sum\limits_{i=1}^{m-1}n_i^2=\left(r+\sum\limits_{i=1}^{m-1}n_i^2\right)^2$.
\end{proof}

If we specialize these general results to rings of the from $\mathcal{R}=\mathbb{Z}/ \mathbb{Z}n=:\mathbb{Z}_n$ 
then we have the following corollaries:

\begin{corollary}
  $$
    \mathcal{I}(\mathbb{Z}_n,1)=n\text{ and }\mathcal{I}(\mathbb{Z}_2,m)=2^m.
  $$
\end{corollary}

\begin{corollary}
  For coprime integers $a$ and $b$ we have 
  $\mathcal{I}(\mathbb{Z}_{ab},m)=\mathcal{I}(\mathbb{Z}_a,m)\cdot\mathcal{I}(\mathbb{Z}_b,m)$.
\end{corollary}

\begin{corollary}
  For a prime $p>2$ we have 
  $$
    \mathcal{I}(\mathbb{Z}_{p^r},m)\ge p^r\cdot p^{m-1\left\lfloor\frac{r}{2}\right\rfloor}.
  $$
\end{corollary}

To be able to do some algebraic calculations later on we denote the set of invertible elements of $\mathcal{R}$ by $\mathcal{R}^*$ and derive a ring $\mathcal{R}'$ from the module $\mathcal{R}^2$.

\begin{definition}
  $$
    \mathcal{R}':=\mathcal{R}[x]/(x^2+1).
  $$
\end{definition}

With $\ii$ being a root of $x^2+1$ we have the following bijection
$$
  \varrho:\mathcal{R}^2\rightarrow\mathcal{R}',\,\,(a,b)\mapsto a+b\ii.
$$
The big advantage of the ring $\mathcal{R}'$ is that we naturally have an addition and multiplication. The construction of the ring is somewhat a reverse engineering of the connection between Paley graphs of square order and integral point sets over the affine plane $\mathbb{F}_q^2$ for $q\equiv 3\mod 4$. With the similar construction of the complex numbers in mind we define:

\begin{definition}
  $$
    \overline{a+b\ii}=a-b\ii.
  $$
\end{definition}

\begin{lemma}
  For $p,p_1,p_2\in\mathcal{R}'$ we have 
  \begin{enumerate}
    \item $d^2(p_1,p_2)=(p_1-p_2)\cdot\overline{(p_1-p_2)}$,
    \item $p\overline{p}\in\mathcal{R}$,
    \item $\overline{p_1+p_2}=\overline{p_1}+\overline{p_2}$,
    \item $\overline{p_1\cdot p_2}=\overline{p_1}\cdot\overline{p_2}$, and
    \item $\overline{\overline{p}}=p$.
  \end{enumerate}
\end{lemma}


\section{Automorphism group of the plane $\mathbf{\mathcal{R}^2}$}
\label{sec_automorphism_group}

Since we want to classify maximal integral point sets up to isomorphism we have to define what we consider as an automorphism. 

\begin{definition}
  \label{def_equivalence}
  An automorphism of $\mathcal{R}'$ with respect to $\Delta$ is a bijective mapping $\varphi$ of $\mathcal{R}'$ with
  \begin{enumerate}
    \item $\Delta(a+b\ii,c+d\ii)=\Delta(\varphi(a+b\ii),\varphi(c+d\ii))$ and\\[-4mm]
    \item there exist $a',b',c',d'\in\mathcal{R}$ such that\\
          $\{\varphi(a+b\ii+r(c+d\ii))\mid r\in\mathcal{R}\}=\{a'+b'\ii+r(c'+d'\ii)\mid r\in\mathcal{R}\}$
  \end{enumerate}
  for all $a$, $b$, $c$, $d$ in $\mathcal{R}$.
\end{definition}

In words this definition says that $\varphi$ has to map points to points, lines to lines, and has to preserve the integral distance property. There is a natural similar definition for $\mathcal{R}^2$ instead of $\mathcal{R}'$.

\begin{lemma}
  \label{lemma_equivalence}
  We have the following examples of automorphisms:
  \begin{enumerate}
    \item $\varphi_s(r)=r+s$ for $s\in\mathcal{R'}$,
    \item $\widetilde{\varphi}(a+b\ii)=b+a\ii$,
    \item $\widetilde{\varphi}_y(r)=ry$ for $y\in{\mathcal{R}'}^*$ with $\exists
          r'\in\mathcal{R}^*:y\overline{y}={r'}^2$, and
    \item $\widehat{\varphi}_j(a+b\ii)=a^{p^j}+b^{p^j}\ii$ for $j\in\mathbb{N}$ and $p$ being the characteristic of
          a field $\mathcal{R}$.
  \end{enumerate}
\end{lemma}
\begin{proof} 
  The first two cases are easy to check. For the third case we consider
  \begin{eqnarray*}
    d^2(r_1y,r_2y)&=&(r_1y-r_2y)\cdot\overline{(r_1y-r_2y)},\\
                  &=& (r_1-r_2)\cdot\overline{(r_1-r_2)}y\overline{y},\\
                  &=&d^2(r_1,r_2)\cdot y\overline{y}.
  \end{eqnarray*}
  For the fourth case we have 
  \begin{eqnarray*}
    d^2(\widehat{\varphi}_j(a_1+b_1\ii),\widehat{\varphi}_j(a_2+b_2\ii))&=&
    (a_1^{p^j}-a_2^{p^j})^2+(b_1^{p^j}-b_2^{p^j})^2,\\
    &=&(a_1-a_2)^{p^j\cdot 2}+(b_1-b_2)^{p^j\cdot 2},\\
    &=&\left((a_1-a_2)^2+(b_1-b_2)^2\right)^{p^j},\\
    &=&d^2(a_1+b_1\ii,a_2+b_2\ii)^{p^j}\\
  \end{eqnarray*}	
  Thus integral point sets are mapped onto integral point sets. That lines are mapped onto lines can be
  checked immediately. Since we have requested that $\mathcal{R}$ is a field for the forth case the mappings
  are injective.
\end{proof}

After this general definition of automorphisms we specialize to the case $\mathcal{R}=\mathbb{F}_q$ with $2\nmid q$. As shorthand we use $\square_q:=\square_{\mathbb{F}_q}$. We remark that the case (4) is the set of Frobenius automorphisms of the field $\mathbb{F}_q$ which is a cyclic group of order $r$ for $q=p^r$.

\begin{theorem}
  \label{thm_automorphism_f_q}
  For $q=p^r$, $p\neq 2$, $q\neq 5,9$ the automorphisms of $\mathbb{F}_q'$ with respect to $\Delta$ are completely 
  described in Lemma \ref{lemma_equivalence}.
\end{theorem}

For $q\equiv 3\mod 4$ this is a well known result on the automorphism group of Paley graphs as mentioned in the introduction. If we consider the set of automorphisms from Lemma \ref{lemma_equivalence} in $\mathbb{F}_q^2$ instead of $\mathbb{F}_q'$ then they form a group with its elements being compositions of the following four mappings:
\begin{enumerate}
  \item $\begin{pmatrix}x\\y\end{pmatrix}\mapsto \begin{pmatrix}x\\y\end{pmatrix} +
        \begin{pmatrix}a\\b\end{pmatrix}$ where $a,b\in\mathbb{F}_q$,\\[1mm]
  \item $\begin{pmatrix}x\\y\end{pmatrix}\mapsto\begin{pmatrix}a&b\\-b&a\end{pmatrix}
        \cdot\begin{pmatrix}x\\y\end{pmatrix}$
        where $a,b\in\mathbb{F}_q$, $a^2+b^2\in\square_q\backslash\{0\}$,\\[1mm]
  \item $\begin{pmatrix}x\\y\end{pmatrix}\mapsto
        \begin{pmatrix}0&1\\1&0\end{pmatrix}\cdot\begin{pmatrix}x\\y\end{pmatrix}$, and \\[1mm]
  \item $\begin{pmatrix}x\\y\end{pmatrix}\mapsto \begin{pmatrix}x^p\\y^p\end{pmatrix}$.
\end{enumerate}

In the remaining part of this section we will prove Theorem \ref{thm_automorphism_f_q}. For the sake of completeness we also give the proof for $q\equiv 3\mod 4$.
If we forget about respecting $\Delta$ then the automorphism group of $\mathbb{F}_q^2$ is the well known group $\text{A$\Gamma$L}(2,\mathbb{F}_q)$. It is a semi-direct product of the translation group, the Frobenius group $\text{Aut}(\mathbb{F}_q)$, and $\text{GL}(2,\mathbb{F}_q)$, the group of multiplications with invertible $2\times 2$ matrices over $\mathbb{F}_q$. So if $G'$ is the automorphism group of $\mathbb{F}_q^2$ with respect to $\Delta$ it suffices to determine the group $G:=G'\cap \text{GL}(2,\mathbb{F}_q)$ because every translation and every element in $\text{Aut}(\mathbb{F}_q)$ respects $\Delta$. So all elements of $G$ can be written as $\begin{pmatrix}x\\y\end{pmatrix}\mapsto\begin{pmatrix}x&y\end{pmatrix}\cdot M$ with $M$ being an invertible $2\times 2$-matrix. As a shorthand we say that $M$ is an element of the automorphism group $G$.

\begin{lemma}
  \label{lemma_sum_of_squares}
  If $M=\begin{pmatrix}a&b\\c&d\end{pmatrix}$ is an element of the automorphism group $G$ then we have
  $ad-bc\neq 0$ and $a^2+b^2$, $a^2+c^2$, $b^2+d^2$, $c^2+d^2$ $\in$ $\square_q$.
\end{lemma}
\begin{proof}
  Since $M$ is also an element of $\text{GL}(2,\mathbb{F}_q)$ its determinant does not vanish. By considering the 
  points $(0,0)$ and $(0,y)$ which are at an integral distance we obtain that $b^2+d^2$ must be a square in $\mathbb{F}_q$.
  Similarly we obtain that $a^2+c^2$, $a^2+b^2$, and $c^2+d^2$ must be squares in $\mathbb{F}_q$.
\end{proof}

To go on we need some facts about roots in $\mathbb{F}_q$ and the set of solutions of quadratic equations in $\mathbb{F}_q$.

\begin{definition}
  For $p^r\equiv 1\mod 4$ we denote by $\omega_q$ an element with $\omega_q^2=-1$.
\end{definition}

\begin{lemma}
  \label{lemma_root_of_minus_one}
  For a finite field $\mathbb{F}_q$ with $q=p^r$ and $p\neq 2$ we have $-1\in\square_q$ iff $q\equiv 1\mod 4$, 
  $\omega_q\in\square_q$ iff $q\equiv 1\mod 8$, and $2\in\square_q$ iff $q\equiv \pm 1\mod 8$.
\end{lemma}
\begin{proof}
  The multiplicative group of the units $\mathbb{F}_q^*$ is cyclic of order $q-1$. Elements of order $4$ are exactly 
  those elements $x$ with $x^2=-1$. A similar argument holds for the the fourth roots of $-1$. The last statement is 
  the second Erg\"anzungssatz of the quadratic reciprocity law generalized to $\mathbb{F}_q$. For a proof we may 
  consider the situation in $\mathbb{F}_p$ and adjungate $x$ modulo the ideal $(x^2-2)$.
\end{proof}

\begin{lemma}
  \label{lemma_quadratic_equation}
  For a fix $c\neq 0$ and $2\nmid q$ the equation $a^2+b^2=c^2$  in $\mathbb{F}_q$ has exactly $q+1$ different 
  solutions if $-1\not\in\square_q$ and exactly $q-1$ different solutions if $-1\in\square_q$.
\end{lemma}
\begin{proof}
  If $b=0$ then we have $a=\pm c$. Otherwise
  $$
    a^2+b^2=c^2\quad\Leftrightarrow\quad \frac{a-c}{b}\cdot\frac{a+c}{b}=-1.
  $$
  Here we set $t:=\frac{a+c}{b}\in\mathbb{F}_q^*$ ($t=0$ corresponds to $b=0$). We obtain
  $$
    2\frac{a}{b}=t-t^{-1},\quad 2\frac{c}{b}=t+t^{-1}\neq 0,
  $$
  yielding
  $$
    t^2\neq -1,\quad b=\frac{2c}{t+t^{-1}},\,\text{ and }\,a= c\cdot\frac{t-t^{-1}}{t+t^{-1}}.
  $$
  If $t$ and $t'$ yield an equal $b$ then we have $t'=t^{-1}$. For $t\neq t^{-1}$ we have different values for 
  $a$ in these cases. Summing up the different solutions proves the stated result. 
\end{proof}

\begin{lemma}
  \label{lemma_cyclic}
  In $\mathbb{F}_q'$ the set $C=\{z\in\mathbb{F}_q'\mid z\overline{z}=1\}$ forms a cyclic multiplicative group.
\end{lemma}
\begin{proof}
  If $-1\not\in\square_q$ then $\mathbb{F}_q'$ is a field and thus $C$ must be cyclic. For the case 
  $-1\in\square_q$ we utilize the bijection
  $$
    \rho_q:\mathbb{F}_q^*\rightarrow G,\quad t\mapsto \frac{1+t^2}{2t}+\omega_q\frac{1-t^2}{2t}x.
  $$
  Now we only have to check that the mapping is a group isomorphism, namely
  $$
    \rho_q(i\cdot j)=\rho_q(i)\cdot\rho_q(j).
  $$
\end{proof}

Our next ingredient is a classification of the subgroups of the projective special linear group $\text{PSL}(2,q)$.

\begin{theorem}(Dickson \cite{dickson}) The subgroups of $\text{PSL}(2,p^r)$ are isomorphic to one of the following families of groups:
  \begin{enumerate}
    \item elementary abelian $p$-groups,
    \item cyclic group of order $z$, where $z$ is a divisor of $\frac{p^r\pm 1}{k}$ and $k=\gcd(p^r-1,2)$,
    \item dihedral group of order $2z$, where $z$ is defined as in (2),
    \item alternating group $A_4$ (this can occur only for $p>2$ or when $p=2$ and $r\equiv 0\mod 12$),
    \item symmetric group $S_4$ (this can only occur if $p^{2r}\equiv 1\mod 16$),
    \item alternating group $A_5$ (for $p=5$ or $p^{2r}\equiv 1\mod 5$),
    \item a semidirect product of an elementary abelian group of order $p^m$ with a cyclic group of order $t$, where 
          $t$ is a divisor of $p^m-1$ and of $p^r-1$, or
    \item the group $\text{PSL}(2,p^m)$ for $m$ a divisor of $r$, or the group $\text{PGL}(2,p^m)$ for $2m$ a 
          divisor of $r$.		
  \end{enumerate}
  \label{theorem_dickson}
\end{theorem}

By $Z:={\pm E}$ we denote the center of $\text{SL}(2,q)$, where $E$ is the identity matrix. Our strategy is to consider $H:=(G\cap\text{SL}(2,q))/Z=G\cap\text{PSL}(2,q)$ and to prove $H\simeq H'$ for $q\ge 13$ where $H'$ is the group of those automorphisms of Lemma \ref{lemma_equivalence} which are also elements of $\text{PSL}(2,q)$. For $-1\not\in\square_q$ we set $\tilde{H}:=\left\{\begin{pmatrix}a&b\\-b&a\end{pmatrix}\mid a^2+b^2=1\right\}$ and for $-1\in\square_q$ we set
$\tilde{H}:=\left\{\begin{pmatrix}a&b\\-b&a\end{pmatrix}\mid a^2+b^2=1\right\} \cup \left\{\begin{pmatrix}-b&a\\a&b\end{pmatrix}\mid a^2+b^2=-1\right\}$.

\begin{lemma}
  \label{lemma_tilde_H}
  For $q\equiv 3\mod 4$ we have $\tilde{H}\simeq\mathbb{Z}_{q+1}$ and for $q\equiv 1\mod 4$ we have 
 $\tilde{H}\simeq D_{q-1}$, where $D_{q-1}$ is the dihedral group of order $2(q-1)$.
\end{lemma}
\begin{proof}
  Utilizing Lemma \ref{lemma_quadratic_equation} and checking that both sets are groups we get
  $$
    |\tilde{H}|=\left\{\begin{array}{ll}q+1&\text{if } q\equiv 3\mod 4,\\[1mm]
    2(q-1)&\text{if }q\equiv 1\mod 4.\end{array}\right.
  $$
  In the first case the group is cyclic due to Lemma \ref{lemma_cyclic}. In the second case it contains a cyclic subgroup 
  of order $q-1$. By checking the defining relations of a dihedral group we can conclude $\tilde{H}\simeq D_{q-1}$ 
  for $q\equiv 1\mod 4$.
\end{proof}

Now we define $H':=\tilde{H} / Z$.

\begin{lemma}
  For $q\ge 13$, $q\equiv 3\mod 4$ we have $H'\simeq\mathbb{Z}_{\frac{q+1}{2}}$ and for $q\ge 13$, $q\equiv 1\mod 4$ 
  we have $H'\simeq D_{\frac{q-1}{2}}$.
\end{lemma}
\begin{proof}
  We have $|H'|=\frac{|\tilde{H}|}{2}$. It remains to show that $H'$ is not abelian for $q\equiv 1\mod 4$. 
  Therefore we may consider the sets $\{\pm M_1\}$ and $\{\pm M_2\}$ where $a,b,c,d$ are elements of $\mathbb{F}_q^*$ 
  with $a^2+b^2=1$, $c^2+d^2=-1$ and where 
  $$
    M_1=
    \begin{pmatrix}
       a & b\\
      -b & a
    \end{pmatrix}
    \quad\text{and}\quad M_2=
    \begin{pmatrix}
      -d & c\\
       c & d
    \end{pmatrix}.
  $$
\end{proof}

\begin{lemma}
  \label{lemma_step}
  For $q\ge 13$ we have $H\simeq H'$.
\end{lemma}
\begin{proof}
  Since $H$ is a subgroup of $\text{PSL}(2,q)$ we can utilize Theorem \ref{theorem_dickson}. We run through the subgroups 
  of $\text{PSL}(2,q)$, identify $H'$ and show that $H$ is no of the subgroups of $\text{PSL}(2,q)$ containing $H'$ as 
  a proper subgroup. With the numbering from the theorem we have the following case distinctions. We remark that for
  $q\equiv 1\mod 4$ the group $H'$ is the group of case (3) and for $q\equiv 3\mod 4$ the group $H'$ is the group of 
  case (2)
  \begin{enumerate}
    \item $H$ is not an elementary abelian $p$-group since $|H'|$ is not a $p$-power.
    \item For $q\equiv 1\mod 4$ the order of $H'$ is larger than $\frac{p^r\pm 1}{2}$ and for $q\equiv 3\mod 4$ 
          the characterized group must be $H'$ itself.
    \item For $q\equiv 1\mod 4$ the characterized group must be $H'$ itself due to the order of the groups. 
          For $q\equiv 3\mod 4$ we must have a look at the elements of order $2$ in $\text{PSL}(2,q)$. These are 
          elements $M\cdot Z$ where $M=\begin{pmatrix}a&b\\c&b\end{pmatrix}$ with $ad-bc=1$ and $M^2=E$ 
          or $M^2=-E$. Solving this equation system yields $M=\pm E$ which corresponds to an element of $H'$ and
          $M=\begin{pmatrix}a&b\\-\frac{a^2+1}{b}&-a\end{pmatrix}$ where $a\in\mathbb{F}_q$ and 
          $b\in\mathbb{F}_q^*$. Now we choose a matrix $N=\begin{pmatrix}u&v\\-v&u\end{pmatrix}$ with
          $u^2+v^2=1$ and $u,v\neq 0$. So $N\cdot Z=\{\pm N\}\in H'$ and since $\langle H',N\rangle$
          would be a dihedral group we have the following relation 
          $$
            MZ \cdot NZ\cdot MZ=N^{-1}Z
          $$
          $$
            \Leftrightarrow\quad
            \{\pm M\}\cdot \{\pm N\}\cdot \{\pm M\}=\{\pm N^{-1}\}=\left\{\pm
            \begin{pmatrix}u&-v\\v&u\end{pmatrix}\right\}
          $$
          $$
            \Leftrightarrow\quad \left\{\pm
            \begin{pmatrix}
              \frac{-ab^2v-a^3v-av-bu}{b}&-v(a^2+b^2)\\
              \frac{v(a^2b^2+a^4+2a^2+1)}{b^2} & \frac{-bu+ab^2v+a^3v+av}{b}
            \end{pmatrix}
            \right\}=\left\{\pm\begin{pmatrix}u&-v\\v&u\end{pmatrix}
            \right\}.
          $$
          By comparing the diagonal elements we get $av(a^2+b^2+1)=0$ and $v(b^4-a^4-2a^2-1)=0$. Due to $v\neq 0$ 
          this is equivalent to $a(a^2+b^2+1)=0$ and $(a^2+b^2+1)\cdot(a^2-b^2+1)=0$. Together with 
          $a^2+b^2\in\square_q$ we conclude $a=0$ and $b=\pm 1$. Since these solutions correspond to an element of 
          $H'$ we derive that case (3) is not possible for $q\equiv 3\mod 4$.
    \item If $H'<H\le A_4$ then $H'$ must be contained in a maximal subgroup of $A_4$. Since the order of a 
          maximal subgroup of $A_4$ is at most $4$ and $q\ge 13$ this case can not occur.
    \item Since we have $q\ge 13$ and the maximal subgroups of the $S_4$ are isomorphic to $A_4$, $D_4$, and $S_3$,
          this case can not occur.
    \item The maximal subgroups of $A_5$ are isomorphic to $D_{5}$, $S_3$, and $A_4$. So this case can not occur 
          for $q\ge 13$.
    \item We have that $|H|$ divides $(q-1)\cdot p^m$. Since $\gcd\left(\frac{q+1}{2},(q-1)\cdot p^m\right)\le 2$ 
          and $|H'|$ divides $|H|$, only $q\equiv 1\mod 4$, $|H'|=q-1$, $t=q-1$, and $r|m$ is possible. If  $m\ge 2r$ 
          then $|H|\ge q^2(q-1)> |\text{PSL}(2,q)|=\frac{1}{2}(q^2-1)q$, which is a contradiction. So only $m=r$ is 
          possible and $H$ must be the semidirect product of an abelian group of order $q$ and a cyclic group of 
          order $q-1$. Using Zassenhaus' theorem \cite[I.18.3]{huppert} we can deduce that all subgroups of order 
          $q-1$ of $H$ are conjugates and so isomorphic. Since $H'$ is not abelian (for $q\equiv 1\mod 4$) it is 
          not cyclic and so at the end case (7) of Theorem \ref{theorem_dickson} is impossible.
    \item Clearly $H\not\simeq \text{PSL}(2,q)$. Since $|H'|$ does not divide 
          $|\text{PSL}(2,p^m)|=\frac{(p^{2m}-1)p^m}{2}$ only the second possibility is left. Since $|H'|$ 
          divides $|\text{PGL}(2,p^m)|=(p^{2m}-1)(p^{2m}-p^m)$ we have $2m=r$, $p^m=\sqrt{q}$, and 
          $q\equiv 1\mod 4$. But for $q\ge 13$ we have $D_{\frac{q-1}{2}}\not\le \text{PGL}(2,\sqrt{q})$, 
          see i.e. \cite{cameron}, thus case (8) is also not possible.
  \end{enumerate}
\end{proof}

To finish the proof of the characterization of the automorphisms of $\mathbb{F}_q^2$ with respect to $\Delta$ we need as a last ingredient a result on the number of solutions of an elliptic curve in $\mathbb{F}_q$.

\begin{theorem}(Hasse, i.e. \cite{02168577})
  \label{thm_hasse}
  Let $f$ be a polynomial of degree $3$ in $\mathbb{F}_q$ without repeated factors then we have for the number 
  $N$ of different solutions of $f(t)=s^2$ in $\mathbb{F}_q^2$ the inequality $|N-q-1|\le 2\sqrt{q}$.
\end{theorem}

\noindent
\textbf{Proof of Theorem \ref{thm_automorphism_f_q}.}
  For the cases $q=3,7,11$ we utilize a computer to check that there are no other automorphisms. So we can assume
  $q\ge 13$.

  If $M\in G$ is an automorphism for $q\equiv 3 \mod 4$ then there exists an element $x\in\mathbb{F}_q^*$ so that 
  either $x\cdot M$ or $x\cdot M\cdot\begin{pmatrix}0&1\\1&0\end{pmatrix}$ has determinant $1$. Thus 
  with the help of Lemma \ref{lemma_step} and Lemma \ref{lemma_equivalence} the theorem is proven for 
  $q\equiv 3\mod 4$. With the same argument we can show that for $q\equiv 1\mod 4$ any possible further automorphism 
  which is not contained in the list of Lemma \ref{lemma_equivalence} must have a determinant which is a non-square 
  in $\mathbb{F}_q$. Let $M=\begin{pmatrix}a&b\\c&d\end{pmatrix}$ be an element of $G$ 
  with $\det(M)=ad-bc\not\in\square_q$. So $M^2=\begin{pmatrix}a^2+bc&b(a+d)\\
  c(a+d)&bc+d^2\end{pmatrix}$ is also an element of $G$. Since we have $\det(M^2)=\det(M)^2\in\square_q$ we 
  have $a^2+bc=bc+d^2$, $b(a+d)=-c(a+d)$ or $a^2+bc=-(bc+d^2)$, $b(a+d)=c(a+d)$ due to Lemma \ref{lemma_step}. This 
  leads to the four cases
  \begin{enumerate}
    \item $a=d$, $b=-c$,
    \item $a=d=0$,
    \item $a=-d$, and
    \item $b=c$, $a^2+d^2=-2b^2$.
  \end{enumerate}
  Now we consider the derived matrix $M':=M\cdot\begin{pmatrix}0&1\\1&0\end{pmatrix}
  =\begin{pmatrix}b&a\\d&c\end{pmatrix}$ with $\det(M')\not\in\square_q$ which must be also 
  an automorphism. So each of the matrices $M$ and $M'$ must be one of the four cases. From this we can conclude 
  some equations and derive a contradiction for each possibility. Here we assume that the number of the case of
  $M'$ is at least the number of the case of $M$.

  \begin{enumerate}
    \item $M$ as in (1): With the help of Lemma \ref{lemma_sum_of_squares} we get
          $\det(M)=a^2+b^2\in\square_q$, which is a contradiction.
    \item $M$ as in (2): Since $\det(M)\not\in\square_q$ the only possibility for $M'$ is case (4). Thus we have 
          $b^2+c^2=0$ $\Leftrightarrow$ $b=\pm \omega_q c$, where we can assume $c=1$ and $b=\omega_q$ without loss 
          of generality. Since $\det(M')$ must be a non-square in $\mathbb{F}_q$ we have $q\equiv 5\mod 8$.
          If we apply $M'$ onto the points $(0,0)$ and $(1,1)$ then we can conclude that $2$ must be a square in
          $\mathbb{F}_q$, which is not the case if $q\equiv 5\mod 8$.
    \item $M$ as in (3): Due to $\det(M)\not\in\square_q$ the matrix $M'$ must be in case (4). So we have 
          $a=d=0$, a situation already treated in case (2).
    \item $M$ as in (4): Thus also $M'$ has to be in case (4). Here we have $a=d$, $b=c$, $2a^2=-2b^2$. Without 
          loss of generality we can assume $a=1$ and $b=\omega_q$. Due to $\det(M)=2\not\in\square_q$ we have
          $q\equiv 5\mod 8$. For two elements $x,y\in\mathbb{F}_q$ with $x^2+y^2$ being a square we have that 
          also $\tilde{M}:=\begin{pmatrix}1&\omega_q\\\omega_q&1\end{pmatrix}
          \cdot\begin{pmatrix}x&y\\-y&x\end{pmatrix}=\begin{pmatrix}x-\omega_q
          x\,\,&\,\,x\omega_q+y\\x\omega_q-y\,\,&\,\,x+y\omega_q\end{pmatrix}$ is an automorphism. Thus 
          with Lemma \ref{lemma_sum_of_squares} we get that $(x\omega_q+y)^2+(x+y\omega)^2=2^2xy\omega_q$ must 
          be a square in $\mathbb{F}_q$ for all possible values $x,y\neq 0$. So for $q\equiv 5\mod 8$ for all possible 
          $x,y$ the product $xy\neq 0$ must be a non-square. We specialize to $x^2+y^2=1^2$ and so can get with the 
          help of Lemma \ref{lemma_quadratic_equation} that $x=\frac{2}{t+t^-1}$ and $y=\frac{t-t^-1}{t+t^-1}$ with 
          $t^2\neq -1$, $t\neq 0$. If we require $t^4\neq 1$ instead of $t^2\neq -1$ we get $x,y\neq 0$. Thus
          $xy=\frac{2(t-t^{-1})}{(t+t^{-1})^2}$ must be a non-square for all $t\in\mathbb{F}_q^*$ with $t^4\neq 1$. 
          Since $2$ is a non-square we have that $t-t^{-1}$ and so also $t^3-t=t(t-1)(t+1)$ must be a square for all
          $t\in\mathbb{F}_q^*$ with $t^4\neq 1$. By checking the five excluded values we see that $f(t):=t(t-1)(t+1)$ 
          must be a square for all $t\in\mathbb{F}_q$. So $f(t)=s^2$ has exactly $N:=2q-3$ solutions in 
          $\mathbb{F}_q$. Since $f$ has not repeated factors and degree $3$ we can apply Theorem \ref{thm_hasse} to 
          get a contradiction to $q\ge 13$.
  \end{enumerate}
\hfill{$\square$}

  %
  %

\begin{lemma}	
  \label{lemma_transitive}
  For two points $p_1\neq p_2\in\mathbb{F}_q'$ at integral distance there exists an isomorphism $\varphi$ with 
  either $\varphi(p_1)=0$, $\varphi(p_2)=1$ or $\varphi(p_1)=0$, $\varphi(p_2)=1+\omega_q\ii$.
\end{lemma}
\begin{proof}
  Without loss of generality we assume $p_1=0$. Since the points $p_1$ and $p_2$ are at integral distance there exists 
  an element $r\in\mathbb{F}_q$ with $p_2\overline{p_2}=r^2$ and since $p_2\neq p_1$ we have $p_2\in{\mathbb{F}_q'}^*$.  
  If $p_2\overline{p_2}\neq 0$ we choose $\cdot p_2^{-1}$ as the isomorphism $\varphi$. Otherwise we have 
  $p_2=a+b\ii$ with $a^2+b^2=0$ where $a,b\neq 0$. Thus $\left(\frac{b}{a}\right)^2=-1$ and $\varphi=\cdot a^{-1}$.
\end{proof}

We remark that Lemma \ref{lemma_transitive} can be sharpened a bit. For three pairwise different non-collinear points $p_1,p_2,p_3\in\mathbb{F}_q'$ with pairwise integral distances there exists an isomorphism $\varphi$ with 
$\{0,1\}\subset\{\varphi(p_1),\varphi(p_2),\varphi(p_3)\}$.

Via a computer calculation we can determine the automorphism groups of the missing cases $q=5,9$.

\begin{lemma}
  For $q=5$ the group $G\le \text{GL}(2,\mathbb{F}_5)$ is given by 
  $$
    \left\{M=\begin{pmatrix}a&b\\\pm b&\pm a\end{pmatrix}\mid
    a,b\in\mathbb{F}_5,\,a^2+b^2\in\square_5,\,\det(M)\neq 0\right\} 
  $$
  where the two signs can be chosen independently.
\end{lemma}


\begin{lemma}
  For $q=9$ the group $G\le \text{GL}(2,\mathbb{F}_9)$ is given by 
  $$
    \left\langle\left\{M=\begin{pmatrix}a&b\\\pm b&\pm a\end{pmatrix}\mid
    a,b\in\mathbb{F}_9,\,a^2+b^2\in\square_9,\,\det(M)\neq 0\right\},
    \begin{pmatrix}1&0\\0&y^2\end{pmatrix}
    \right\rangle 
  $$
  where the two signs can be chosen independently and where $y$ is a primitive root in $\mathbb{F}_9^*$.
\end{lemma}

For $q=5$ there are exactly $32$ such matrices and for $q=9$ there are exactly $192$ such matrices. For $q=5,9$ Lemma \ref{lemma_transitive} can be sharpend. Here the automorphism group acts transitively on the pairs of points with integral 
distance, as for $q\equiv 3\mod 4$.

We would like to remark that also for $q\equiv 3\mod 4$ the automorphism group of $\mathbb{F}_q^2$ with respect to $\Delta$ is isomorphic to the automorphism group of the quadrance graph over $\mathbb{F}_q^2$. This can easily be verified be going over the proof of Theorem \ref{thm_automorphism_f_q} again and by checking the small cases using a computer.


\section{Maximal integral point sets in the plane $\mathbb{F}_q^2$}
\label{sec_maximum_integral_point_sets}

\noindent
Very nice rings are those which are integral domains. These are in the case of finite commutative rings exactly 
the finite fields $\mathbb{F}_q$ where $q=p^r$ is a prime power. So far we only 
have the lower bound $\mathcal{I}(\mathbb{F}_q,2)\ge q$. In this section we will prove $\mathcal{I}(\mathbb{F}_q,2)=q$ 
for $q>2$. In the case of $\mathbb{F}_p$ we will even classify the maximum integral point sets up to isomorphism. One way to prove $\mathcal{I}(\mathbb{F}_q,2)=q$ for $2\nmid q$ is to consider the graph $\mathcal{G}_q$ with the elements of $\mathbb{F}_q$ as its vertices and pairs of points at integral distance as edges. For $q\equiv 3\mod 4$ the graph $\mathcal{G}_q$ is isomorphic to the Paley graph of order $q^2$. From \cite{0561.12009} we know that in this case a maximum clique of $\mathcal{G}_q$ has size $q$ and is isomorphic to a line. Also for $q\equiv 1\mod 4$ the graph $\mathcal{G}_q$ is a strongly regular graph. So we can apply a result from \cite{0423.05019,0466.05026} on cliques of strongly regular graphs. It turns out that a maximum clique has size $q$ and that every clique $\mathcal{C}$ of size $q$ is \textit{regular}, in the sense of \cite{0423.05019,0466.05026}, this means in our special case that every point not in $\mathcal{C}$ is adjacent to $\frac{q+1}{2}$ points in $\mathcal{C}$. To start with our classification of maximum integral point sets over $\mathbb{F}_q$ we need the concept of directions.

\begin{definition}
  For a point $p=a+b\ii\in\mathbb{F}_q'$ the quotient $\frac{b}{a}\in\mathbb{F}_q\cup\{\infty\}$ is called
  the \textbf{direction} of $p$. For two points $p_1=a_1+b_1\ii$, $p_2=a_2+b_2\ii$ the direction is defined as
  $\frac{b_1-b_2}{a_1-a_2}\in\mathbb{F}_q\cup\{\infty\}$. We call an direction $d$ \textbf{integral} if two points
  $p_1$, $p_2$ with direction $d$ have an integral distance.
\end{definition}

Point sets of cardinality $q$ in $\mathbb{F}_q^2$ with at most $\frac{q+3}{2}$ directions are more or less completely classified:

\begin{theorem}{(Ball, Blokhuis, Brouwer, Storme, Sz\H{o}nyi, \cite{0945.51002})} 
  \label{thm_directions}
  Let $f:\mathbb{F}_q\rightarrow\mathbb{F}_q$, 
  where $q=p^n$, $p$ prime, $f(0)=0$. Let $N=|D_f|$, where $D_f$ is the set of directions determined by the function $f$.
  Let $e$ (with $0\le e\le n$) be the largest integer such that each line with slope in $D_f$ meets the graph of $f$ in a 
  multiple of $p^e$ points. Then we have the following:
  \begin{enumerate}
    \item $e=0$ and $\frac{q+3}{2}\le N\le q+1$,
    \item $e=1$, $p=2$, and $\frac{q+5}{3}\le N\le q-1$,
    \item $p^e>2$, $e|n$, and $\frac{q}{p^e}+1\le N\le\frac{q-1}{p^e-1}$,
    \item $e=n$ and $N=1$.
  \end{enumerate}
  Moreover, if $p^e>3$ or ( $p^e=3$ and $N=\frac{q}{3}+1$), then $f$ is a linear map on $\mathbb{F}_q$ viewed as 
  a vector space over $\mathbb{F}_{p^e}$. (All possibilities for $N$ can be determined in principle.)
\end{theorem}

Here a function $f:\mathbb{F}_q\rightarrow\mathbb{F}_q$ determines a point set $\Ps=\{(x,f(x))\mid x\in\mathbb{F}_q\}$ of cardinality $q$. In the case $N=1$ the point set is a line. In the case $e=0$ and $N=\frac{q+3}{2}$ then $\Ps$ is affine equivalent to the point set corresponding to $x\mapsto x^{\frac{q+1}{2}}$.

We remark that affine equivalence is a bit more than our equivalence because we have to respect $\Delta$. The next thing to prove is that integral point sets can not determine too many directions.

\begin{lemma}
  \label{lemma_different_directions}
  For $2\nmid q$ an integral point set over $\mathbb{F}_q^2$ determines at most $\frac{q+3}{2}$ different directions 
  if $-1\in\square_q$ and at most $\frac{q+1}{2}$ different directions if $-1\not\in\square_q$.
\end{lemma}
\begin{proof}
  We consider the points $p=a+b\ii$ at integral distance to $0$. Thus there exists an element $c'\in\mathbb{F}_q$ 
  with $a^2+b^2=c'^2$. In the case $a=0$ we obtain the direction $\infty$. Otherwise we set $d:=\frac{b}{a}$ and 
  $c:=\frac{c'}{a}$, yielding $1=c^2-d^2=(c-d)(c+d)$, where $d$ is the direction of the point. Now we set 
  $c+d=:t\in\mathbb{F}_q^*$ yielding $c=\frac{t+t^{-1}}{2}$, $d=\frac{t-t^{-1}}{2}$. The two values $t$ and $-t^{-1}$ 
  produce an equal direction. Since $t=-t^{-1}$ $\Leftrightarrow$ $t^2=-1$ we get the desired bounds.
\end{proof}

We need a further lemma on the number of points on a line in a non collinear integral point set:

\begin{lemma}
  \label{lemma_points_on_line}
  If $2\nmid q$ and $\Ps$ is a non collinear integral point set over $\mathbb{F}_q^2$, 
  then each line $l$ contains at most $\frac{q-1}{2}$ points for $-1\notin\square_q$ and 
  at most $\frac{q+1}{2}$ points for $-1\in\square_q$.
\end{lemma}
\begin{proof}
  If $l$ is a line with an integral pair of points on it, then its slope is an integral direction. 
  Now we consider the intersections of lines with integral directions containing a point $p\notin
  l$, with $l$.
\end{proof}

We remark that there would be only $\frac{q-1}{2}$ integral directions for $q\equiv 1\mod 4$ if we would not consider $0$ as a square as for quadrance graphs. In this case there could be at most $\frac{q-3}{2}$ points on $l$ for $q\equiv 1\mod 4$ in Lemma \ref{lemma_points_on_line}.


To completely classify maximum integral point sets over $\mathbb{F}_q'$ we need the point set $\mathcal{P}_q:=(1\pm\omega_q\ii)\square_q$.

\begin{lemma}
  $\mathcal{P}_q$ is an integral point set of cardinality $q$.
\end{lemma}
\begin{proof}
  \begin{eqnarray*}
    \quad\quad\quad\quad\quad
    d^2(r_1^2+r_1^2\omega_q\ii,r_2^2+r_2^2\omega_q\ii)&=&0^2,\\
    d^2(r_1^2+r_1^2\omega_q\ii,r_2^2-r_2^2\omega_q\ii)&=&(2\omega_qr_1r_2)^2,\\
    d^2(r_1^2-r_1^2\omega_q\ii,r_2^2-r_2^2\omega_q\ii)&=&0^2.\\
  \end{eqnarray*}
\end{proof}

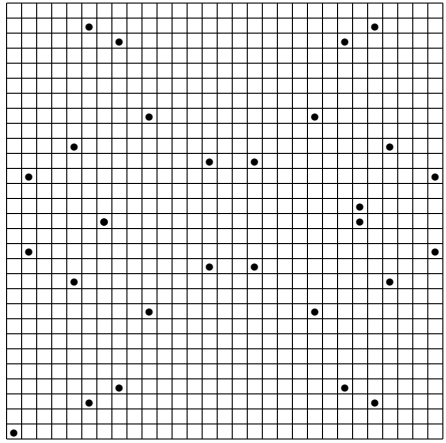
\begin{figure}[ht]
  \begin{center}
    \setlength{\unitlength}{0.2cm}
    \begin{picture}(29,29)
      \multiput(0,0)(0,1){30}{\line(1,0){29}}
      \multiput(0,0)(1,0){30}{\line(0,1){29}}
      \put(0.15,0.15){\tiny{$\bullet$}}
      \put(1.15,12.15){\tiny{$\bullet$}}
      \put(4.15,19.15){\tiny{$\bullet$}}
      \put(5.15,2.15){\tiny{$\bullet$}}
      \put(6.15,14.15){\tiny{$\bullet$}}
      \put(7.15,26.15){\tiny{$\bullet$}}
      \put(9.15,21.15){\tiny{$\bullet$}}
      \put(13.15,11.15){\tiny{$\bullet$}}
      \put(16.15,18.15){\tiny{$\bullet$}}
      \put(20.15,8.15){\tiny{$\bullet$}}
      \put(22.15,3.15){\tiny{$\bullet$}}
      \put(23.15,15.15){\tiny{$\bullet$}}
      \put(24.15,27.15){\tiny{$\bullet$}}
      \put(25.15,10.15){\tiny{$\bullet$}}
      \put(28.15,17.15){\tiny{$\bullet$}}
      \put(1.15,17.15){\tiny{$\bullet$}}
      \put(4.15,10.15){\tiny{$\bullet$}}
      \put(5.15,27.15){\tiny{$\bullet$}}
      \put(6.15,14.15){\tiny{$\bullet$}}
      \put(7.15,3.15){\tiny{$\bullet$}}
      \put(9.15,8.15){\tiny{$\bullet$}}
      \put(13.15,18.15){\tiny{$\bullet$}}
      \put(16.15,11.15){\tiny{$\bullet$}}
      \put(20.15,21.15){\tiny{$\bullet$}}
      \put(22.15,26.15){\tiny{$\bullet$}}
      \put(23.15,14.15){\tiny{$\bullet$}}
      \put(24.15,2.15){\tiny{$\bullet$}}
      \put(25.15,19.15){\tiny{$\bullet$}}
      \put(28.15,12.15){\tiny{$\bullet$}}
    \end{picture}
    \caption{The maximum integral point set $\mathcal{P}_{29}$.}
    \label{fig_the_cross}
  \end{center}
\end{figure}

In Figure \ref{fig_the_cross} we have depicted $\mathcal{P}_{29}$ as an example. By construction the points of $\mathcal{P}_q$ are located on the two lines $(1,\omega_q)\cdot\mathbb{F}_q$ and $(1,-\omega_q)\cdot\mathbb{F}_q$ which intersect in $(0,0)$ with an \textit{angle} of $90$ degree, but this fact seems not that obvious by looking at Figure \ref{fig_the_cross}. We remark that this construction of $\mathcal{P}_q$ works in any commutative ring $\mathcal{R}$ where $-1\in\square_{\mathcal{R}}$ and that none of these point sets corresponds to a quadrance graph. If we apply this construction on $\mathcal{R}=\mathbb{Z}_{p^r}$ we obtain an integral point set of cardinality $\phi(p^r)+1=(p-1)\cdot p^{r-1}+1$, where $\phi$ is the Euler-function defined by $\phi(n)=|\mathbb{Z}_n^*|$.

\begin{lemma}
  For $2\mid r$ the point set $\Ps:=\{(a,b)\mid a,b\in\mathbb{F}_{\sqrt{q}}\}$ is an 
  integral point set.
\end{lemma}
\begin{proof}
  We have $\mathbb{F}_{\sqrt{q}}\subset\square_q$.
\end{proof}

We remark that for $\sqrt{q}\equiv 1\mod 4$ also the point set $\Ps:=\{(a,\omega_q b)\mid a,b\in\mathbb{F}_{\sqrt{q}}\}$ is integral.

We say that an integral point set is maximal if we can not add a further point without destroying the property \textit{integral point set}. All given examples of integral point sets of size $q$ are maximal. This could be proved be applying results on cliques of strongly regular graphs or in the following way.

\begin{lemma}
  The lines $1\cdot\mathbb{F}_q$ and $(1+\omega_q\ii)\cdot\mathbb{F}_q$ are maximal.
\end{lemma}
\begin{proof}
  We apply Lemma \ref{lemma_points_on_line}.
\end{proof}

\begin{lemma}
  The integral point set $\Ps=(1\pm\omega_q\ii)\cdot \square_q$ is maximal.
\end{lemma}
\begin{proof}
  \label{lemma_max_2}
  Let us assume there is a further point $(a+b\ii)\not\in\Ps$ with $a,b\in\mathbb{F}_q$ such that $\Ps\cup\{(a+b\ii)\}$ 
  is also an integral point set. We know that $(a+b\ii)$ can not lie on one of the lines 
  $(1+\omega_q \ii)\cdot\mathbb{F}_q$ or $(1-\omega_q \ii)\cdot\mathbb{F}_q$. Thus $a^2+b^2\neq 0$. The points of $\Ps$ 
  are given by $(1+\omega_q\ii)r_1^2$ and $(1-\omega_q\ii)r_2^2$ for arbitrary $r_1,r_2\in\mathbb{F}_q$. We define 
  functions $f_1,f_2:\mathbb{F}_q\rightarrow\mathbb{F}_q$ via
  \begin{eqnarray*}
    f_1(r_1)=(a-r_1^2)^2+(b-r_1^2\omega_q)^2=a^2+b^2-2r_1^2(a+b\omega_q),\\
    f_2(r_2)=(a-r_2^2)^2+(b+r_2^2\omega_q)^2=a^2+b^2-2r_2^2(a-b\omega_q).\\
  \end{eqnarray*}
  Since these are exactly the squared distances of the points of $\Ps$ to the point $(a+b\ii)$ we have
  $\text{Bi}(f_1),\text{Bi}(f_2)\subseteq \square_q$. Using a counting argument we have $\text{Bi}(f_1),
  \text{Bi}(f_2)= \square_q$. The term $-2(a+b\omega_q)$ is a fix number. Let us assume that it is a square. 
  Then for each square $r^2$ and $c=a^2+b^2\neq 0$ the difference $r^2-c$ must be a square. But the equation 
  $r^2-c=h^2$ has $\frac{q+1}{2}<q$ solutions for $r$, which is a contradiction. Thus $-2(a+b\omega_q)$ and
  $-2(a-b\omega_q)$ are non-squares. But $r^2-c\not\in \square_q$ has $\frac{q-1}{2}$ solutions, thus we
  have a contradiction 
\end{proof}

\begin{theorem}
  \label{main_thm}
  For $q=p^r>9$ with $p\neq 2$, $r=1$ or $q\equiv 3\mod 4$ an integral point set of
  cardinality $q$ is isomorphic to one of the stated examples.
\end{theorem}
\begin{proof}
  We consider a point set $\Ps$ of $\mathbb{F}_{q}$ of cardinality $q$ with at most 
  $\frac{q+3}{2}$ directions and utilize Theorem \ref{thm_directions}. If $e=r$ and 
  $N=1$ then $\Ps$ is a line. If $e=1$ then $\Ps$ is affine equivalent to 
  $X:=\{(x,x^{\frac{q+1}{2}})\mid x\in\mathbb{F}_q\}$. This is only possible for $q\equiv 1 \mod 4$.
  The set $X$ consists of two orthogonal lines. Since there are only two types of non-isomorphic 
  integral lines in $\mathbb{F}_q^2$ and each point $p$ not on a line $l$ is at integral distance to
  $\frac{q+1}{2}$ points on $l$ we have two unique candidates of integral point sets of this type.
  One is given by $(1\pm\omega_q\ii)\cdot\square_q$. For the other possibility we may assume that
  $(0,0),(1,0)\in\Ps$. Thus $(0,\pm\omega_q)\in\Ps$, $(-1,0),(\pm\omega_q,0),(0,\pm 1)\in\Ps$. So
  $\Ps$ must be symmetric in the following sense: There exists a set $S\subset\mathbb{F}_q^*$ such
  that $\Ps=(0,0)\cup\{(0,a),(a,0)\mid a\in S\}$. The elements $s$ of $S$ must fullfill 
  $s\in\mathbb{F}_q^*$, $s^2+1\in\square_q$ and $s^2-1\in\square_q$. Each condition alone has only $\frac{q-1}{2}$
  solutions. Fulfilling both conditions, meaning $|S|=\frac{q-1}{2}$ is
  possible only for $q\le 9$. For $q=5,9$ there are such examples.
  For $q\equiv 3\mod 4$ we refer to \cite{0561.12009}.
\end{proof}

We remark that there may be further examples of integral point sets of cardinality $q$ for $q=p^r\equiv 1\mod 4$ and $r>1$. Those examples would correspond to case (3) of Theorem \ref{thm_directions}.

\begin{theorem}
  \label{thm_cardinality}
  For $q=p^r$ with $p\neq 2$ we have $\mathcal{I}(\mathbb{F}_{q},2)=q$.
\end{theorem}
\begin{proof}
  Let $\Ps$ be an arbitrary integral point set of cardinality $q$. Now we show that $\Ps$ is maximal. If 
  we assume that there is another integral point set $\Ps'$ with $\Ps\subset\Ps'$ and $|\Ps'|=q+1$ 
  then we can delete a point of $\Ps'$ in such a way that we obtain an integral point set $\Ps''$ 
  with $e=1$ in the notation of \ref{thm_directions}. Thus $\Ps''\simeq(1\pm\omega_q\ii)\cdot\square_q$.
  Since $\Ps''$ is maximal due to Lemma \ref{lemma_max_2} we have a contradiction.
\end{proof}



\section{Maximal integral point sets in the plane $\mathbb{Z}_n^2$}
\label{sec_maximal_integral_point_sets_modular}

\noindent
Due to Theorem \ref{thm_decompose} for the determination of $\mathcal{I}(\mathbb{Z}_n,2)$ we only need to consider the cases $n=p^r$.

\begin{lemma}
  \label{lemma_2}
  $$
    \mathcal{I}(\mathbb{Z}_{p^{r+1}},2)\le p^2\cdot \mathcal{I}(\mathbb{Z}_{p^r},2).
  $$
\end{lemma}
\begin{proof}
  We consider the natural ring epimorphism $\nu:\mathbb{Z}_{p^{r+1}}\rightarrow\mathbb{Z}_{p^r}$. If $\Ps$ is an
  integral point set in $\mathbb{Z}_{p^{r+1}}^2$ then $\nu(\Ps)$ is an integral point set in $\mathbb{Z}_{p^r}^2$.
\end{proof}

For $p\ge 3$ we have the following examples of integral point sets in $\mathbb{Z}_{p^r}^2$ with big cardinality (with some abuse of notation in the third case).
\begin{eqnarray*}
  &&\left\{\left(i,j\cdot p^{\left\lceil\frac{r}{2}\right\rceil}\right)\mid i,j\in\mathbb{Z}_{p^r}\right\},\\
  &&\left\{\left(i,i\omega_{\mathbb{Z}_{p^r}}+j\cdot p^{\left\lceil\frac{r}{2}\right\rceil}\right)
  \mid i,j\in\mathbb{Z}_{p^r}\right\},\text{ and}\\
  && (1,\pm\omega_{\mathbb{Z}_{p}})\cdot \square_{\mathbb{Z}_{p}}+\left\{(p\cdot
  a,p\cdot b)\mid a,b\in\mathbb{Z}_{p^r}\right\}\text{ for }r=2.
\end{eqnarray*}
Each of these examples has cardinality $p^r\cdot p^{\left\lfloor\frac{r}{2}\right\rfloor}$.

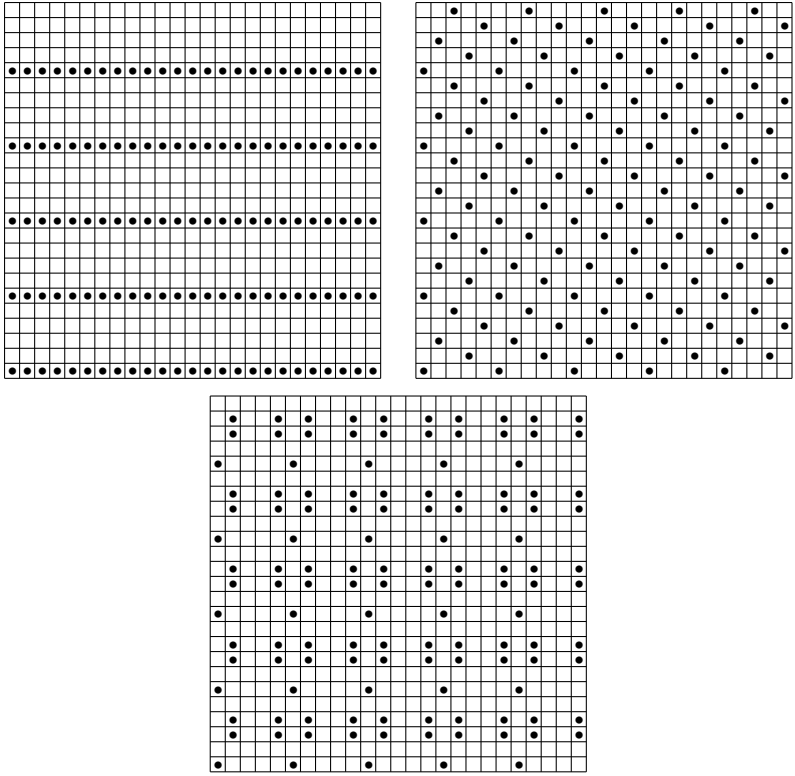
\begin{figure}[ht]
  \begin{center}
    \setlength{\unitlength}{0.2cm}
    \begin{picture}(25,25)
      \multiput(0,0)(0,1){26}{\line(1,0){25}}
      \multiput(0,0)(1,0){26}{\line(0,1){25}}
      \multiput(0,0)(1,0){25}{\put(0.2,0.2){\tiny{$\bullet$}}}
      \multiput(0,5)(1,0){25}{\put(0.2,0.2){\tiny{$\bullet$}}}
      \multiput(0,10)(1,0){25}{\put(0.2,0.2){\tiny{$\bullet$}}}
      \multiput(0,15)(1,0){25}{\put(0.2,0.2){\tiny{$\bullet$}}}
      \multiput(0,20)(1,0){25}{\put(0.2,0.2){\tiny{$\bullet$}}}
    \end{picture}
    \quad
    \begin{picture}(25,25)
      \multiput(0,0)(0,1){26}{\line(1,0){25}}
      \multiput(0,0)(1,0){26}{\line(0,1){25}}
      \multiput(0,0)(0,5){5}{\put(0.2,0.2){\tiny{$\bullet$}}}
      \multiput(0,0)(0,5){5}{\put(1.2,2.2){\tiny{$\bullet$}}}
      \multiput(0,0)(0,5){5}{\put(2.2,4.2){\tiny{$\bullet$}}}
      \multiput(0,0)(0,5){5}{\put(3.2,1.2){\tiny{$\bullet$}}}
      \multiput(0,0)(0,5){5}{\put(4.2,3.2){\tiny{$\bullet$}}}
      \multiput(0,0)(0,5){5}{\put(5.2,0.2){\tiny{$\bullet$}}}
      \multiput(0,0)(0,5){5}{\put(6.2,2.2){\tiny{$\bullet$}}}
      \multiput(0,0)(0,5){5}{\put(7.2,4.2){\tiny{$\bullet$}}}
      \multiput(0,0)(0,5){5}{\put(8.2,1.2){\tiny{$\bullet$}}}
      \multiput(0,0)(0,5){5}{\put(9.2,3.2){\tiny{$\bullet$}}}
      \multiput(0,0)(0,5){5}{\put(10.2,0.2){\tiny{$\bullet$}}}
      \multiput(0,0)(0,5){5}{\put(11.2,2.2){\tiny{$\bullet$}}}
      \multiput(0,0)(0,5){5}{\put(12.2,4.2){\tiny{$\bullet$}}}
      \multiput(0,0)(0,5){5}{\put(13.2,1.2){\tiny{$\bullet$}}}
      \multiput(0,0)(0,5){5}{\put(14.2,3.2){\tiny{$\bullet$}}}
      \multiput(0,0)(0,5){5}{\put(15.2,0.2){\tiny{$\bullet$}}}
      \multiput(0,0)(0,5){5}{\put(16.2,2.2){\tiny{$\bullet$}}}
      \multiput(0,0)(0,5){5}{\put(17.2,4.2){\tiny{$\bullet$}}}
      \multiput(0,0)(0,5){5}{\put(18.2,1.2){\tiny{$\bullet$}}}
      \multiput(0,0)(0,5){5}{\put(19.2,3.2){\tiny{$\bullet$}}}
      \multiput(0,0)(0,5){5}{\put(20.2,0.2){\tiny{$\bullet$}}}
      \multiput(0,0)(0,5){5}{\put(21.2,2.2){\tiny{$\bullet$}}}
      \multiput(0,0)(0,5){5}{\put(22.2,4.2){\tiny{$\bullet$}}}
      \multiput(0,0)(0,5){5}{\put(23.2,1.2){\tiny{$\bullet$}}}
      \multiput(0,0)(0,5){5}{\put(24.2,3.2){\tiny{$\bullet$}}}
    \end{picture}\\[2mm] 
    \begin{picture}(25,25)
      \multiput(0,0)(0,1){26}{\line(1,0){25}}
      \multiput(0,0)(1,0){26}{\line(0,1){25}}
      \multiput(0,0)(5,0){5}{\multiput(0,0)(0,5){5}{\put(0.2,0.2){\tiny{$\bullet$}}}}
      \multiput(0,0)(5,0){5}{\multiput(0,0)(0,5){5}{\put(1.2,2.2){\tiny{$\bullet$}}}}
      \multiput(0,0)(5,0){5}{\multiput(0,0)(0,5){5}{\put(1.2,3.2){\tiny{$\bullet$}}}}
      \multiput(0,0)(5,0){5}{\multiput(0,0)(0,5){5}{\put(4.2,2.2){\tiny{$\bullet$}}}}
      \multiput(0,0)(5,0){5}{\multiput(0,0)(0,5){5}{\put(4.2,3.2){\tiny{$\bullet$}}}}
    \end{picture}
    \caption{Three maximal integral point sets over $\mathbb{Z}_{25}^2$ of cardinality $125$.}
    \label{fig_maximal_z_p_2}
  \end{center}
\end{figure}

\begin{conjecture}
  \label{conj_2}
  The above list is the complete list of maximum integral point sets in $\mathbb{Z}_{p^r}^2$ up to isomorphism.
\end{conjecture}

So far we do not even know the automorphism group of $\mathbb{Z}_n^2$ with respect to $\Delta$. But with Definition \ref{def_equivalence} Conjecture \ref{conj_2} is well defined. Using Lemma \ref{lemma_equivalence} we know at least a subgroup of the automorphism group. If there are any further automorphisms is an open question which has to be analyzed in the future.

\begin{theorem}
  For $p\ge 3$ we have $\mathcal{I}(\mathbb{Z}_{p^2},2)=p^3$ and the above list of extremal examples is complete.
\end{theorem}
\begin{proof}
  With $\mathcal{I}(\mathbb{Z}_{p},2)=p$, Lemma \ref{lemma_2} and the examples we get
  $\mathcal{I}(\mathbb{Z}_{p^2},2)=p^3$.
  Let $\Ps$ be a maximum integral point set in $\mathbb{Z}_{p^2}$. By $S$ denote the lower left $p\times p$-square of 
  $\mathbb{Z}_{p^2}$
  $$
    S:=\{(i,j)+\mathbb{Z}_{p^2}^2\mid 0\le i,j\le p-1,\,i,j\in\mathbb{Z}\}.
  $$
  Using Theorem \ref{thm_cardinality} and Lemma \ref{lemma_2} we can deduce that for each
  $(u,v)\in\mathbb{Z}_{p^2}^2$ we have
  $$
    |\Ps\cap\left((u,v)+S\right)|\le p.
  $$
  Since we can tile $\mathbb{Z}_{p^2}$ with $p^2$ such sets (including $S+(u,v)$) equality must hold.
  After a transformation we can assume that $\Ps\cap S$ equals one of the three
  following possibilities
  \begin{enumerate}
    \item $\{(i,0)\mid 0\le i\le p-1\}$,
    \item $\{(i,\omega_{\mathbb{Z}_p} i)\mid 0\le i\le p-1\}$, or
    \item $(1,\pm\omega_{\mathbb{Z}_p})\cdot \square_{\mathbb{Z}_p}$.
  \end{enumerate}
  In the first case we consider $\Ps\cap (S+(1,0))$. With Lemma \ref{lemma_points_on_line} we get $(p,0)\in\Ps$ 
  and iteratively we get $(i,0)\in\Ps$ for all $i\in\mathbb{Z}_{p^2}$. Now we consider $\Ps\cap (S+(0,1))$ and 
  conclude $\Ps=\left\{\left(i,j\cdot p\right)\mid i,j\in\mathbb{Z}_{p^2}\right\}$. With the same argument we can 
  derive $\Ps=\left\{\left(i,i\omega_{\mathbb{Z}_p}+j\cdot p\right)\mid i,j\in\mathbb{Z}_{p^2}\right\}$ in the second 
  case and $\Ps= (1,\pm\omega_{\mathbb{Z}_{p}})\cdot \square_{\mathbb{Z}_{p}}+\left\{(p\cdot
  a,p\cdot b)\mid a,b\in\mathbb{Z}_{p^r}\right\}$ in the third case.
\end{proof}


\section{Maximal integral point sets with no three collinear points}
\label{section_arcs}

\noindent
In this and the next section we study the interplay between the integrality condition for a point set and further common restrictions for lines and circles.

\begin{definition}
  A set of $r$ points $(u_i,v_i)\in\mathcal{R}^2$ is said to be \textbf{collinear} if there are
  $a,b,t_1,t_2,w_i\in\mathcal{R}$ with
  $$
    a+w_it_1=u_i\quad\text{and}\quad b+w_it_2=v_i.
  $$
\end{definition}

There is an easy necessary criterion to decide whether three points are collinear.

\begin{lemma}
  \label{lemma_collinear}
  If three points $(u_1,v_1)$, $(u_2,v_2)$, and $(u_3,v_3)\in\mathcal{R}^2$ are collinear then it holds
  $$
    \left|\begin{pmatrix}
      u_1&v_1&1\\
      u_2&v_2&1\\
      u_3&v_3&1\\
    \end{pmatrix}\right|=0.
  $$
\end{lemma}

If $\mathcal{R}$ is an integral domain the above criterion is also sufficient. The proof is easy and left to the reader.

\begin{definition}
  By $\overline{\mathcal{I}}(\mathcal{R},2)$ we denote the maximum cardinality of an integral point set with no three 
  collinear points.
\end{definition}

\begin{lemma}
  $$
    \overline{\mathcal{I}}(\mathcal{R},2)\le 2\cdot |\mathcal{R}|.
  $$
\end{lemma}
\begin{proof}
  We ignore the integrality condition and consider the lines $l_i=\{(i,r)\mid r\in\mathcal{R}\}$ for
  all $i\in\mathcal{R}$.
\end{proof}

\begin{lemma}
  \label{lemma_upper_bound_circle}
  If $-1\in\square_q$ we have $\overline{\mathcal{I}}(\mathbb{F}_q,2)\le \frac{q+3}{2}$ and for $-1\not\in
  \square_q$ we have $\overline{\mathcal{I}}(\mathbb{F}_q,2)\le \frac{q+1}{2}$.
\end{lemma}
\begin{proof}
  Let $\Ps$ be an integral point set over $\mathbb{F}_q$ without a collinear triple. We choose a point $p\in\Ps$. 
  The directions of $p$ to the other points $p'$ of $\Ps$ are pairwise different. Since there are at most $\frac{q+3}{2}$ 
  or $\frac{q+1}{2}$ different directions in an integral point set over $\mathbb{F}_q$ 
  (Lemma \ref{lemma_different_directions}), we obtain $|\Ps|\le\frac{q+5}{2}$ for $-1\in\square_q$ and
  $|\Ps|\le\frac{q+3}{2}$ for $-1\not\in\square_q$. Suppose that this upper bound is achieved. So all points must have
  exactly one \textit{neighbor} in direction $0$ and one in direction $\infty$. Thus $|\Ps|$ must be even in this case,
  which is a contradiction due to Lemma \ref{lemma_root_of_minus_one}.
\end{proof}

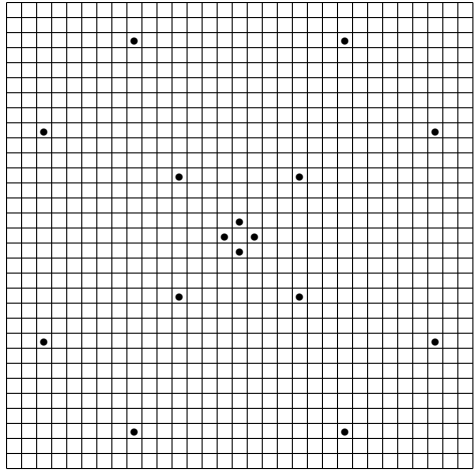
\begin{figure}[ht]
  \begin{center}
    \setlength{\unitlength}{0.2cm}
    \begin{picture}(31,31)
      \multiput(0,0)(0,1){32}{\line(1,0){31}}
      \multiput(0,0)(1,0){32}{\line(0,1){31}}
      \put(2.15,8.15){\tiny{$\bullet$}}
      \put(2.15,22.15){\tiny{$\bullet$}}
      \put(8.15,2.15){\tiny{$\bullet$}}
      \put(8.15,28.15){\tiny{$\bullet$}}
      \put(11.15,11.15){\tiny{$\bullet$}}
      \put(11.15,19.15){\tiny{$\bullet$}}
      \put(14.15,15.15){\tiny{$\bullet$}}
      \put(15.15,14.15){\tiny{$\bullet$}}
      \put(15.15,16.15){\tiny{$\bullet$}}
      \put(16.15,15.15){\tiny{$\bullet$}}
      \put(19.15,11.15){\tiny{$\bullet$}}
      \put(19.15,19.15){\tiny{$\bullet$}}
      \put(22.15,2.15){\tiny{$\bullet$}}
      \put(22.15,28.15){\tiny{$\bullet$}}
      \put(28.15,8.15){\tiny{$\bullet$}}
      \put(28.15,22.15){\tiny{$\bullet$}}
    \end{picture}
    \caption{Integral point set corresponding to the construction from Lemma \ref{lemma_circle_construction} for $q=31$.}
    \label{fig_semi_general_position}
  \end{center}
\end{figure}

Using an element $z\in\mathcal{R}'$ with $z\overline{z}=1$ we can describe a good construction for lower bounds. Actually this equation describes something like a circle with radius one. An example for $q=31$ is depicted in Figure \ref{fig_semi_general_position}.

\begin{lemma}
  \label{lemma_circle_construction}
  For $z\in\mathcal{R}'$ with $z\overline{z}=1$ the set $\Ps=\{z^{2i}\mid i\in\mathbb{N}\}$ is an integral point set.
\end{lemma}
\begin{proof}
  With $c:=a-b$ we have
  \begin{eqnarray*}
    d(z^{2a},z^{2b})&=& (z^{2a}-z^{2b})\cdot\overline{(z^{2a}-z^{2b})}=(z^{2c}-1)\cdot\overline{z^{2c}-1}\\
    &=&2-z^{2c}2-\overline{z^{2c}}=(\underset{\in\mathcal{R}}{\underbrace{z^c\ii-\overline{z^c}\ii}})^2
  \end{eqnarray*}
\end{proof}

We remark that the set $\Ps'=\{z^{2i+1}\mid i\in\mathbb{N}\}$ is an isomorphic integral point set. The set of solutions of 
$z\overline{z}=1$ forms a cyclic multiplicative group $G$ due to Lemma \ref{lemma_cyclic}. From Lemma \ref{lemma_quadratic_equation} we know that $G$ has size $q+1$ for $-1\not\in\square_q$ and size $q-1$ if $-1\in\square_q$. So by Lemma \ref{lemma_circle_construction} we get a construction of an integral point set in $\mathbb{F}_q$ which is near the upper bound of Lemma \ref{lemma_upper_bound_circle}. We only have to prove that our construction does not produce three collinear points in $\mathbb{F}_q$.

\begin{lemma}
  For $\mathcal{R}=\mathbb{F}_q$ with $2\nmid q$ the point set from Lemma \ref{lemma_circle_construction} 
  contains no collinear triple.
\end{lemma}
\begin{proof}
  We assume that we have three pairwise different points $p_1,p_2,p_3$ in $\mathcal{R}'$ which are collinear. 
  So there exist $a,b,c,d,t_1,t_2,$ and $t_3$ in $\mathcal{R}$ fullfilling
  \begin{eqnarray*}
    p_1&=&a+bt_1+(c+dt_1)\ii,\\
    p_2&=&a+bt_2+(c+dt_2)\ii,\\
    p_3&=&a+bt_3+(c+dt_3)\ii,\\
  \end{eqnarray*}
  and $t_i\neq t_j$ for $i\neq j$. Since $p_i\overline{p_i}=1$ we have
  \begin{eqnarray*}
    a^2+2abt_1+b^2t_1^2+c^2+2cdt_1+d^2t_1^2 &=& 1,\\
    a^2+2abt_2+b^2t_2^2+c^2+2cdt_2+d^2t_2^2 &=& 1,\\
    a^2+2abt_3+b^2t_3^2+c^2+2cdt_3+d^2t_3^2 &=& 1.\\
  \end{eqnarray*}
  Subtracting the first two and the last two equations yields
  \begin{eqnarray*}
    2ab(t_1-t_2)+b^2(t_1-t_2)(t_1+t_2)+2cd(t_1-t_2)+d^2(t_1-t_2)(t_1+t_2) &=&0,\\
    2ab(t_2-t_3)+b^2(t_2-t_3)(t_2+t_3)+2cd(t_2-t_3)+d^2(t_2-t_3)(t_2+t_3) &=&0.
  \end{eqnarray*}
  Because $t_1\neq t_2$, $t_2\neq t_3$ and $\mathcal{R}$ is an integral domain we obtain
  \begin{eqnarray*}
    2ab+b^2(t_1+t_2)+2cd+d^2(t_1+t_2) &=&0,\\
    2ab+b^2(t_2+t_3)+2cd+d^2(t_2+t_3) &=&0.
  \end{eqnarray*}
  Another subtraction yields
  \begin{eqnarray*}
    &&b^2(t_1-t_3)+d^2(t_1-t_3)=0\quad
    \Rightarrow\quad b^2+d^2=0.
  \end{eqnarray*}
  Inserting yields
  $$
    2ab+2cd=0\quad\Leftrightarrow\quad 2ab=-2cd
  $$
  and
  $$
    a^2+c^2=1.
  $$
  Thus
  $$
    4a^2b^2=4c^2d^2\quad\Leftarrow\quad (a^2+c^2)4b^2=0\quad\Leftrightarrow\quad b=0.
  $$
  In the same way we obtain $d=0$ and so $p_1=p_2=p_3$, which is a contradiction.
\end{proof}

\begin{corollary}
  For $-1\not\in\square_q$ we have $\overline{\mathcal{I}}(\mathbb{F}_q,2)=\frac{q+1}{2}$ and for $-1\in
  \square_q$ we have $\frac{q-1}{2}\le\overline{\mathcal{I}}(\mathbb{F}_q,2)\le\frac{q+3}{2}$.
\end{corollary}

\begin{conjecture}
  \label{conj_3}
  For $-1\in\square_q$ we have $\overline{\mathcal{I}}(\mathbb{F}_q,2)=\frac{q-1}{2}$.
\end{conjecture}

We remark that Conjecture \ref{conj_3} would be true for quadrance graphs. Following the proof of Lemma \ref{lemma_upper_bound_circle} we would get $\frac{q-1}{2}$ as an upper bound for $q\equiv 1\mod 4$. Since $z^c-\overline{z^c}=0$ would imply $2c=q-1$ the construction from Lemma \ref{lemma_circle_construction} does not contain a pair of points with squared distance $0$.


\section{Integral point sets in general position}
\label{sec_general_position}

\noindent
Our best construction for integral point sets where no three points are collinear consists of points on a \textit{circle}. So it is interesting to study integral point sets where additionally no $4$ points are allowed to be situated on a \textit{circle}.

\begin{definition}
  \label{def_circle}
  Points $p_i=(x_i,y_i)$ in $\mathcal{R}^2$ are said to be situated on a circle if there exist $a,b,r\in\mathcal{R}$ 
  with $(x_i-a)^2+(y_i-b)^2=r$ for all $i$.
\end{definition}

We have the following condition:

\begin{lemma}
  \label{lemma_check_on_circle}
  Four distinct points $p_i=(x_i,y_i)$ in $\mathbb{F}_q^2$ which contain no collinear triple
  are situated on a circle if and only if
  $$
    \left|\begin{pmatrix}
      x_1&y_1&x_1^2+y_1^2&1\\
      x_2&y_2&x_2^2+y_2^2&1\\
      x_3&y_3&x_3^2+y_3^2&1\\
      x_4&y_4&x_4^2+y_4^2&1
    \end{pmatrix}\right|=0.
  $$
\end{lemma}
\begin{proof}
  If there exist $a,b,r\in\mathbb{F}_q$  with $(x_i-a)^2+(y_i-b)^2=r$ for all $1\le i\le 4$ then the determinant 
  clearly vanishes since $r=(x_i-a)^2+(y_i-b)^2=(x_i^2+y_i^2)-2a\cdot x_i-2b\cdot y_i+(a^2+b^2)$. For the other direction
  we consider the unique circle $\mathcal{C}$ through the points $(x_1,y_1)$, $(x_2,y_2)$, $(x_3,y_3)$ described by 
  the parameters $a,b,r\in\mathbb{F}_q$. With the same idea as before we get
  $$
    \left|\begin{pmatrix}
      x_1&y_1&0&1\\
      x_2&y_2&0&1\\
      x_3&y_3&0&1\\
      x_4&y_4&(x_4-a)^2+(y_4-b)^2-r&1
    \end{pmatrix}\right|=0.
  $$
    If $(x_4,y_4)$ is not on the circle $\mathcal{C}$ then we can develop the determinant after the third
    column and obtain
  $$
    \left|\begin{pmatrix}
      x_1&y_1&1\\
      x_2&y_2&1\\
      x_3&y_3&1
    \end{pmatrix}\right|=0.
  $$
  which is a contradiction to the fact that $(x_1,y_1)$, $(x_2,y_2)$, and $(x_3,y_3)$ are not collinear, 
  see Lemma \ref{lemma_collinear}.
\end{proof}

We remark that for arbitrary commutative rings $\mathcal{R}$ the determinant criterion from Lemma \ref{lemma_check_on_circle} is a necessary condition. 

\begin{definition}
  \label{def_general_position}
  By $\dot{\mathcal{I}}(\mathcal{R},2)$ we denote the maximum cardinality of an integral point set in 
  $\mathcal{R}^2$ which is in general position, this means that it contains no collinear triple and no 
  four points on a circle.
\end{definition}

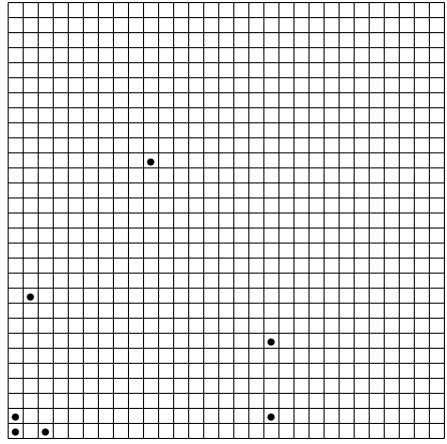
\begin{figure}[ht]
  \begin{center}
    \setlength{\unitlength}{0.2cm}
    \begin{picture}(29,29)
      \multiput(0,0)(0,1){30}{\line(1,0){29}}
      \multiput(0,0)(1,0){30}{\line(0,1){29}}
      \put(0.15,0.15){\tiny{$\bullet$}}
      \put(0.15,1.15){\tiny{$\bullet$}}
      \put(2.15,0.15){\tiny{$\bullet$}}
      \put(9.15,18.15){\tiny{$\bullet$}}
      \put(1.15,9.15){\tiny{$\bullet$}}
      \put(17.15,1.15){\tiny{$\bullet$}}
      \put(17.15,6.15){\tiny{$\bullet$}}
    \end{picture}
    \caption{A maximum integral point set in general position over $\mathbb{F}_{29}^2$.}
    \label{fig_general_position}
  \end{center}
\end{figure}

As a shorthand for the conditions of Definition \ref{def_general_position} we also say that the points are in general position. An example of seven points over $\mathbb{F}_{29}^2$ in general position which pairwise integral distances is depicted in Figure \ref{fig_general_position}. As trivial upper bound we have $\dot{\mathcal{I}}(\mathcal{R},2)\le\overline{\mathcal{I}}(\mathcal{R},2)$. By applying the automorphisms of $\mathbb{F}_q^2$ with respect to $\Delta$ we see that they conserve circles.

\begin{table}[ht]
  \begin{center}
    \begin{tabular}{|r|r||r|r||r|r||r|r||r|r|}
      \hline $n$ & $\IC$ & $n$ & $\IC$ & $n$ & $\IC$ & $n$ &    $\IC$ & $n$ &    $\IC$ \\ 
      \hline   2 &     4 &  17 &     5 &  41 &     9 &  67 &  9 &  97 & 11 \\ 
      \hline   3 &     2 &  19 &     5 &  43 &     8 &  71 & 11 & 101 & 13 \\ 
      \hline   5 &     4 &  23 &     5 &  47 &     7 &  73 & 10 & 103 & 11 \\ 
      \hline   7 &     3 &  29 &     7 &  53 &     9 &  79 & 11 & 107 & 11 \\ 
      \hline  11 &     4 &  31 &     6 &  59 &     9 &  83 & 11 & 109 & 12 \\ 
      \hline  13 &     5 &  37 &     7 &  61 &    10 &  89 & 11 & 113 & 12 \\ 
      \hline 
    \end{tabular} 
    \caption{Values of $\dot{\mathcal{I}}(\mathbb{F}_p,2)=\dot{\mathcal{I}}(\mathbb{Z}_p,2)$ for small primes $p$.}
    \label{table_n4oc}
  \end{center}
\end{table}

Via an exhaustive combinatorial search we have determined $\dot{\mathcal{I}}(\mathbb{F}_p,2)$ for small values of $p$, see Table \ref{table_n4oc}. Since it is a non-trivial task to determine these numbers exactly, at least for $p\ge 100$, we give an outline of our used algorithm.

\begin{algo}{\textbf{(Generation of integral point sets in general position over }$\mathbf{\mathbb{F}_q}$)}\\
  \label{algo_main}
  \textit{Input:} $q$\\
  \textit{Output:} Integral point sets $\Ps\subset\mathbb{F}_q$ in general position\\
  {\bf begin}\\
  \Tsm $\Ps=[(0,0),(0,1)]$\\
  \Tsm $blocked[(0,0)]=blocked[(0,1)]=true$\\
  \Tsm{\bf loop over }$d\in\mathbb{F}_q\,\,\,\mathcal{L}_d=[]$ {\bf end}\\
  \Tsm{\bf loop over }$x\in\mathbb{F}_q^2\backslash\{(0,0),(0,1)\}$\\
  \Tsm\Tsm$\text{blocked}[x]=false$\\
  \Tsm\Tsm{\bf if }$\Delta((0,0),x)=0$ {\bf or }$\Delta((0,1),x)=0$ {\bf then }$\text{blocked}[x]=true$ {\bf end}\\
  \Tsm\Tsm{\bf if }$collinear((0,0),(0,1),x)$ {\bf then } $\text{blocked}[x]=true$ {\bf end}\\
  \Tsm\Tsm{\bf if }$\text{blocked}[x]=true$ {\bf then } $\mathcal{L}_{get\_direction(x)}.append(x)$ {\bf end}\\
  \Tsm{\bf end}\\
  \Tsm $add\_point(\Ps,0)$\\
  {\bf end}
\end{algo}

So far almost nothing is done. We restrict our search to integral point sets $\Ps$ of cardinality at least $3$. So we may assume that $\Ps$ contains the points $(0,0)$ and $(0,1)$. For each $x\in\mathbb{F}_q^2$ the variable $blocked[x]$ says whether $x$ can be appended to $\Ps$ without destroying the property integral point set or general position. The lists $\mathcal{L}_d$ cluster the points of $\mathbb{F}_q^2$ according to their direction. The fact that $\Ps$ can contain besides $(0,0)$ and $(0,1)$ at most one member from each $\mathcal{L}_d$ can be used to prune the search tree if one searches only for integral point sets with maximum cardinality.

\begin{algo}{\textbf{($\mathbf{add\_point}$)}}\\
  \textit{Input:} Lower bound $l$ on the direction and an integral point set $\Ps$\\
  \textit{Output:} Integral point sets $\Ps\subset\mathbb{F}_q$ in general position\\
  {\bf begin}\\
  \Tsm{\bf loop over } $d\in\mathbb{F}_q$ {\bf with } $d\ge l$\\
  \Tsm\Tsm{\bf loop over }$x\in\mathcal{L}_d$ {\bf with } $blocked[x]=false$\\
  \Tsm\Tsm\Tsm{\bf if }$canon\_check(\Ps,x)=true$ {\bf then}\\
  \Tsm\Tsm\Tsm\Tsm$\mathcal{P}.append(x)$\\
  \Tsm\Tsm\Tsm\Tsm block all $y\in\mathbb{F}_q^2$ where $\Delta(y,x)=0$ {\bf or }$collinear(p_1,x,y)=true$\\
  \Tsm\Tsm\Tsm\Tsm{\bf or } $on\_circle(p_1,p_2,x,y)=true$ for $p_1,p_2\in\Ps$\\
  \Tsm\Tsm\Tsm\Tsm output $\Ps$\\
  \Tsm\Tsm\Tsm\Tsm $add\_point(\Ps,d+1)$\\
  \Tsm\Tsm\Tsm\Tsm unblock\\
  \Tsm\Tsm\Tsm\Tsm$\mathcal{P}.remove(x)$\\
  \Tsm\Tsm\Tsm{\bf end}\\
  \Tsm\Tsm{\bf end}\\
  \Tsm{\bf end}\\
  {\bf end}
\end{algo}

The subroutine $add\_point$ simply adds another point to the point set $\Ps$ and maintains the set of further candidates for adding a further point. Some lookahead is possible to implement. Since the automorphism group of $\mathbb{F}_q^2$ with respect to $\Delta$ is very large we would obtain lots of isomorphic integral point sets if we do without isomorphism pruning. With the framework of orderly generation, see i.e. \cite{winner}, it is possible to write a subroutine $canon\_check$ that let our algorithm output a complete list of pairwise non-isomorphic integral point sets in general position. For our purpose it suffices to have a subroutine $canon\_check$ that rejects the majority of isomorphic copies but as a return has a good performance. Let $m:\mathbb{F}_q^2\rightarrow\mathbb{F}_q^2$, $(x,y)\mapsto(-x,y)$ the automorphism that mirrors at the $y$-axis and let $\preceq$ be a total ordering on the points of $\mathbb{F}_q^2$ if $u\prec v$ for $direction(u)<direction(v)$. For the latter comparison we use an arbitrary but fix total ordering of $\mathbb{F}_q$, where $0$ is the smallest element and which is also used for the looping over $\mathbb{F}_q$. By $\Ps[2]$ we denote the third point of a list $\Ps$.

\begin{algo}{\textbf{($\mathbf{canon\_check}$)}}\\
  \textit{Input:} An integral point set $\Ps$\\
  \textit{Output:} Returns false if $\Ps$ should be rejected due to isomorphism pruning\\
  {\bf begin}\\
  \Tsm{\bf loop over} some disjoint triples $(u,v,w)\in\Ps\times\Ps\times\Ps$ {\bf with } $\delta^2(u,v)\neq 0$\\
  \Tsm\Tsm determine an automorphism $\alpha$ with $\alpha(u)=(0,0)$ and $\alpha(v)=(0,1)$\\
  \Tsm\Tsm{\bf if } $\alpha(w)\prec \Ps[2]$ {\bf or } $m(\alpha(w))\prec \Ps[2]$ {\bf then return } $false$ {\bf end}\\
  \Tsm{\bf end}\\
  \Tsm{\bf return } $true$\\
  {\bf end}
\end{algo}

For further examples we refer to \cite{own_hp} where we list the coordinates of one extremal example for $p\le 113$.

A formal proof of the correctness of the proposed algorithm is not difficult but a bit technical and so left to the reader. We remark that there are several non-isomorphic integral point sets in general position which achieve the upper bound $\dot{\mathcal{I}}(\mathbb{Z}_n,2)$. So far we have no insight in their structure or in the asymptotic behavior of. $\dot{\mathcal{I}}(\mathbb{Z}_n,2)$. It seems that we have $\dot{\mathcal{I}}(\mathbb{Z}_p,2)\ge 7$ for all sufficiently large primes $p$. This is interesting because the question whether $\dot{\mathcal{I}}(\mathbb{Z},2)\ge 7$ is unsolved 
so far. In other words, there is no known $7_2$-cluster \cite{UPIN}. This is a set of seven points in the plane, no three points on a line, no four points on a circle, where the coordinates and the pairwise distances are integral.

\begin{conjecture}
  For each $l$ there is a $p'$ so that for all $p\ge p'$ we have $\dot{\mathcal{I}}(\mathbb{Z}_p,2)\ge l$. 
\end{conjecture}

%

\section{Conclusion and outlook}

\noindent
In this paper we have considered sets of points $\Ps$ in the affine plane $\text{AG}(2,q)$ with pairwise integral distances. We have presented several connections to other discrete structures and problems. Some questions concerning maximum cardinalities and complete classifications of extremal examples remain open. Clearly similar questions could be asked in $\text{AG}(3,q)$ or higher dimensional spaces.

\section*{Acknowledgment}
\noindent
I am thankful to Aart Blokhuis, Stancho Diemiev, Michael Kiermaier, Harald Meyer, and Ivo Radloff whose comments were very helpful during the preparation of this article.

\bibliography{integral_over_finite_fields}
\bibdata{integral_over_finite_fields}
\bibliographystyle{plain}

\end{document}